\documentclass{article}

\usepackage{latexsym,amsfonts,amssymb,proof}
\usepackage{fancyvrb}
\usepackage{txfonts}
\usepackage{epsf}

\newcommand{\ket}[1]{|{#1}\rangle}
\newcommand{\av}{``}

\newtheorem{definition}{Definition}[section]
\newtheorem{theorem}[definition]{Theorem}
\newtheorem{lemma}[definition]{Lemma}

\newtheorem{proposition}[definition]{Proposition}

\title {QUANTUM STATES AS VIRTUAL SINGLETONS:\\
CONVERTING DUALITY INTO SYMMETRY}
\author{Giulia Battilotti\\
Department of Mathematics\\
University of Padova - Italy\\}
\date{}

\begin{document}
\maketitle
\begin{abstract}
In a  predicative framework from basic logic,  defined for a model of quantum parallelism by sequents, we characterize a class of first order domains, termed {\em virtual singletons}, which allows a generalization of the notion of duality, termed {\em symmetry}.  Although consistent
with the classical notion of duality, symmetry creates an environment where negation has fixed points, for which the direction of logical consequence is irrelevant. Symmetry can model Bell's states. So, despite its nonsense in a traditional logical setting, symmetry can hide the peculiar advantage for the treatment of information, that is proper of quantum mechanics.
\end{abstract}
\section*{Introduction}
The problem of modelling quantum mechanics by logic is among the main open  topics in the search for quantum structures.  It has been acquiring more and more relevance since the proposal
of quantum computers, that is of exploiting quantum processes for computation and communication processes. Hence the problem gathers different aspects in the foundations of mathematics, of physics, and now of computation. 
 
Then, in the recent field of quantum computational logics, many different approaches are present.
Here, we remind two lines of research which we feel closer to our motivations. 
Some developments are founded on the logical analysis of the algebraic model given by Hilbert spaces. We quote the first survey paper \cite{DCGL03}, moreover we remind that such a choice,  sided  by other different approaches, considering different spaces for quantum mechanics, such as Foch spaces,  or different algebraic structures such as MV-algebras, has permitted the logical formalization of different aspects of the theory to this community of research. 
Considering instead a non-algebraic approach to the problem of quantum computation,  allowing a direct comparison with logical systems, in a computational framework, we quote in particular \cite{BS}, for the embedding of quantum entanglement in a logical system.

The idea that has motivated our research in quantum computational logics is that a logical system for quantum computation need to include a logical representation of quantum entanglement. This is not a harmless requirement, in logic, since it is becoming more and more apparent that the traditional input-output functional way to conceive an algorithm is not suited to model quantum computational processes, that are based on quantum entanglement. For a detailed reading of quantum algorithms in this sense we quote \cite{Ca}, and his reference \cite{DE} for an experiment enlightening the physical motivation of non-sequentiality. Moreover, we quote \cite{BV}, which performs an analysis of open quantum system dynamics and concludes that a quantum algorithm is not the sequence of its temporal parts, due to the entangling gates, that behave as \av black boxes". As is well known, the functional view of algorithms has, in logic, a clear translation: logical implication.  We remind that the semantics of intuitionistic implication is based on the notion of function.  Hence, we find it very intriguing to overcome the traditional functional view of computation in logic itself, discussing the definition of connectives. Discussing in which terms logical connectives are defined, what is lost and what is gained when one abandons the usual setting of logic, which gives room to the definition of logical implication, should offer a good opportunity to face the problem of the meaning of computation.

So, we search the different approach to information that is offered by quantum mechanics, directly in terms of assertions and logical equations defining connectives. This is possible in basic logic \cite{SBF}, a sequent calculus platform for several extensional logics, including some quantum logics. 
In such a framework, we have proposed a characterization of quantum states in predicative logic \cite{Ba}, \cite{Ba2}. 

Here, we introduce and discuss  the notion of \av virtual singleton", namely a set which can act as a singleton under particular
assumptions, interpreting such notion in our model. We see that virtual singletons determine a symmetric rather than dual setting in logic, that means we can achieve a list of 
self-dual literals, that are fixed points for negation and have rather a phase duality. In our model, we can see that 
phase duality is naturally induced by the duality already present in logic. Simply, we are not aware of it, since it occurs prior to measurement. 
 Virtual singletons are insensitive to the direction of
logical consequence and allow a different symmetric connective, obtained as a generalized quantifier, whose computation is \av parallel" and excludes a context-free computation. It permits an interpretation of Bell's states in our model. This constitutes an alternative to logical implication.  

A very intriguing open problem is to understand how the implication could then arise, when 
the intrinsic randomness of QM is dropped in favour of determinism. We think that it should arise when the \av functional view" is recovered, due to the fact that the usual meaning of first order variables can be established, when virtual singletons disappear. 

Virtual singletons are more similar to infinite sets than to finite, countable, sets, and their logical setting offers a holistic rather than compositional treatment of information. Such a feature is present in quantum mechanics, and also in our human thinking, since, for example, it is a feature of the human language. This problem  has already been considered in  quantum computational semantics, for example in \cite{DCGL2}. As for our approach, we have realized that the logical setting derived in our model is surprisingly close to \av bi-logic", introduced by the psychoanalyst I. Matte Blanco in the '70s \cite{MB}. For, in bi-logic, negation and implication are meaningless, and information is \av infinite". These connections, not analyzed in the present paper, will be object of future work, and, we believe, they are a strong motivation for this kind of studies, beyond the modelling of computation.

\section{Symmetry and duality in basic logic}
Basic logic  is a common platform to study sequent calculi for extensional logics. Its sequent calculus, and its extensions, have the following symmetry theorem (\cite{SBF}):
\begin{theorem}
The sequent $\Gamma\vdash \Delta$ is proved by a proof $\Pi$ of the calculus $\bf{B}$ of basic logic if and only if the symmetric sequent $\Delta^s\vdash \Gamma^s$ is proved by the symmetric  proof $\Pi^s$ of $\bf{B}$. 
The equivalence is preserved for any symmetric extension of $\bf{B}$. Moreover, the equivalence is preserved switching \av right" and \av left" extensions
of  $\bf{B}$. 
\end{theorem}
The theorem can be proved putting any involution $(-)^s$ on literals (including the identity), and putting $(A\circ B)^s\equiv
B^s\circ^s A^s$, where the pairs $(\circ, \circ^s)$ for logical constants are the following: the additive conjunction and disjunction $(\&, \vee))$, the multiplicative disjunction and conjunction $(\ast, \otimes)$\footnote{Here we denote the multiplicative disjunction by $\ast$}, implication and exclusion $(\to, \leftarrow)$. Then one proves the statement by induction on proofs:
\begin{itemize}
\item
It is $A^s\vdash A^s$ for every $A$ (axioms are symmetric)
\item
Inference rules  are in symmetric pairs:
$$
\frac{\Gamma_2\vdash \Delta_2 }{\Gamma_0\vdash \Delta_0\quad \Gamma_1\vdash \Delta_1}\circ R/L \quad \Longleftrightarrow \quad 
\frac{\Delta_2^s\vdash \Gamma_2^s }{\Delta_1^s\vdash \Gamma_1^s\quad \Delta_0^s\vdash \Gamma_0^s}\circ L/R
$$
where $\circ^s L/R$ is the right (resp. left) rule for the connective $\circ^s$ when $\circ R/L$ is the left (resp. right) rule for $\circ$.
\end{itemize}
The involution $(-)^s$ is non trivial on the inductive step, namely logical constants are defined in {\em dual} rather than symmetric pairs: $(\&, \vee)\dots$, whereas it can be the identity on axioms. Then the orientation of the turnstyle is irrelevant for axioms, it is relevant only for inference rules. So the orientation of logical consequence and logical theorems is due to the duality, applied at the inductive steps.

In logic, the symmetry theorem gets full meaning when one considers couples of dual literals: $p, p^\perp$ (Girard's literals) an puts the involution $(-)^\perp$: $(p)^\perp\equiv p^\perp$, $(p^\perp)^\perp\equiv p$. Then the statement of the symmetry theorem is written as follows:
\begin{equation}\label{symdual}
\Gamma\vdash \Delta \quad \Longleftrightarrow \quad \Delta^\perp\vdash \Gamma^\perp
\end{equation}
In this form, the duality is not yet a negation. For, in order to have the usual characterization of negation (Girard's negation in our case), we should have a formulation with contexts:
\begin{equation}\label{extdual}
\Gamma',\Gamma\vdash \Delta, \Delta' \quad \Longleftrightarrow \quad \Gamma', \Delta^\perp\vdash \Gamma^\perp, \Delta'
\end{equation}
{\em for every} way to separate $\Gamma, \Gamma'$ and $\Delta, \Delta'$. 

Such a formulation is not derivable from the symmetry theorem. For, even if rules of basic logic can be extended to liberalize contexts at the left and/or at the right, getting calculi for several extensional logics, including linear and classical logic, the symmetry theorem cannot work with contexts.
For example, the sequent $p\to q, p\vdash q$ is provable in some extensions of basic logic, whereas the sequent $p\to q, q\vdash p$ is nowhere provable!
A proof of the sequent $p\to q, q^\perp\vdash p^\perp$, valid in some extensions of basic logic, is not obtained as the symmetric proof of  $p\to q, p\vdash q$, but by suitable structural rules on the duality $\perp$ (\cite{FS}). In such a way extended duality (\ref{extdual}) is proved too. 

Hence the symmetry theorem characterizes a \av proto-negation", with a context-sensitive behaviour. One has usual negation only extending it to a standard context-free calculus. This could mean, possibly, to kill different potentialities of symmetry. 
Our question is: is there a real \av symmetric" interpretation of the symmetry theorem? In which terms is it in conflict with usual duality?
More practically: what is a symmetric literal? Are there symmetric connectives somewhere? The model we have developed for quantum computation can contribute to give an answer.

The answer seems negative in usual propositional logic. For, as proved in basic logic, the usual propositional constants and their rules can be interpreted as derivable from suitable definitory equations, as in the schema:
$$
\Gamma\vdash A\circ^r B\; \equiv\; \Gamma\vdash A\cdot B\qquad\qquad B\circ^l A\vdash \Delta\;\equiv\; B\cdot A\vdash \Delta
$$
In such a schema, the connectives $\circ^r$ and $\circ^l$ are defined as the result of importing a metalinguistic link, represented by $\cdot$, into the formal
language.
So the same metalinguistic link  can define a couple of dual connectives $(\circ^r, \circ^l)$, at the left and at the right of the turnstyle. 
Then the  solution of definitory equations gives symmetric pairs of rules, for the couple $(\circ^r, \circ^l)$, and hence the symmetry theorem sketched above (see \cite{SBF}).
\subsection{Symmetry in the predicative case}
Such a setting works loose in the predicative case, as we are illustrating. 
The meaning of the quantifiers is given by considering assertions linked by the metalinguistic link {\em forall} \cite{MS}.
One has the assertion \av $\Gamma$ {\em yields} $A(z)$", where $z$ is a variable, $\Gamma$ does not depend on $z$ free (in the following, we adopt the notation $\Gamma(-z_1, z_2\dots)$ to say that the variables $z_1, z_2\dots$ are not free in $\Gamma$),  and the free variable, on a range $D$, is the glue {\em for all} the $A(z)$. Then we write our assertion:

\centerline{{\em forall} $z\in D$, \quad $\Gamma(-z)\vdash A(z)$} 

\noindent where the premise $z\in D$ is at a metalinguistic level. We import it into the language of sequents and write:

\centerline{$\Gamma(-z), z\in D \vdash A(z)$}

\noindent We consider it the primitive assertion generating the quantifier $\forall$. So we put the following definitory equation of $\forall$:
\begin{equation}\label{foralldef}
\Gamma(-z)\vdash (\forall x\in D)A(x) \quad{\mbox{\em if and only if}}\quad \Gamma(-z), z\in D\vdash A(z)
\end{equation}
so that the meaning of $(\forall x\in D)A(x)$ is \av$\forall x(x\in D\to A(x))$".

In order to find the {\em symmetric} logical translation of the link, we keep the metalinguistic \av {\em forall} $z\in D$" and consider a symmetric sequent where $A(z)$ is on the left:
$A(z)\vdash \Delta(-z)$. The assertion  

\centerline{{\em forall} $z\in D$, \quad $A(z)\vdash \Delta(-z)$} 
\noindent  importing $z\in D$ into the sequent, has the form 

\centerline{$A(z)\vdash \Delta(-z), z\notin D$}
\noindent since one has that $A(z)$ {\em  yields} $\Delta(-z)$ {\em unless} $z\notin D$. So a negation is required in order to express the symmetric
meaning of the metalinguistic link {\em forall}. In the above sequent, we  substitute, formally the negated proposition $x\notin D$ by the dual proposition
 $(x\in D)^d$, where $d$ is a hypothetical duality, so that  we put the following symmetric definitory equation for the quantifier $\exists$:
$$
(\exists x\in D)A(x)\vdash \Delta(-z) \quad{\mbox{\em if and only if}}\quad A(z)\vdash \Delta(-z), (z\in D)^d
$$
where the meaning of $(\exists x\in D)A(x)$ is attributed by exclusion: \av $\exists x(A(z)\leftarrow (x\in D)^d)$" (exclusion $\leftarrow$ is the symmetric
connective of implication $\rightarrow$ in basic logic). Notice that usually one has no need of $d$ in order to define an existential quantifier in logic. For,  if $d$ satisfies the equivalence characterizing negation: $\Gamma, x\in D\vdash \Delta$ iff $\Gamma\vdash x\notin D, \Delta$, the sequent $A(z)\vdash \Delta(-z), (z\in D)^d$ is equivalent to $A(z), z\in D\vdash \Delta(-z)$. Then
the usual way to define $(\exists x\in D)A(x)$ follows: $\exists x(x\in D\& A(x))$. In it no reference to the duality $d$ is present. 

One can prove that the following symmetric pairs of rules derive from the solution of our definitory equations:
$$
\frac{\Gamma(-z), z\in D \vdash A(z)}
{\Gamma(-z)\vdash (\forall x\in D)A(x)}\forall f
\qquad
\frac{\Gamma \vdash z\in D \quad A(z)\vdash \Delta}
{\Gamma, (\forall x\in D)A(x)\vdash \Delta}\forall r
$$
$$
\frac{A(z)\vdash \Delta(-z), (z\in D)^d}{(\exists x\in D)A(x)\vdash \Delta(-z)}\exists f 
\qquad
\frac{\Gamma\vdash A(z) \quad (z\in D)^d\vdash \Delta}{\Gamma\vdash (\exists x\in D)A(x), \Delta}\exists r
$$
Sketch of proof: $\forall f$ is one direction of the definitory equation and $\forall r$ is the minimal requirement to re-obtain the reflection axiom: $z\in D, (\forall x\in D)A(x)\vdash A(z)$,
in turn derived putting $(\forall x\in D)A(x)$ for $\Gamma(-z)$ in the definitory equation. In the symmetric equation, the reflection axiom is
$A(z)\vdash (\exists x\in D)A(x), (z\in D)^d$. For more information, see \cite{SBF} and \cite{MS}.
The $\forall$ rules are proper of an intuitionistic calculus, the $\exists$ rules are proper of a dual-intuitionistic calculus. In the following, we shall
feel free to enrich them with additional contexts, up to rules of classical logic.

We assume to have a language with equality. Then we consider the following Leibnitz-style definitory equation for the equality predicate $=$, introduced in basic logic by Maietti \cite{Ma}.
For every way to separate $\Gamma, \Gamma'$ and $\Delta, \Delta'$ it is
$$
\Gamma', \Gamma(t/s), s=t \vdash \Delta(t/s), \Delta' \quad \mbox{iff} \quad\Gamma', \Gamma \vdash \Delta, \Delta'
$$
We also write down the symmetric equivalence, that defines the predicate $\neq$ at the right
$$
\Gamma', \Gamma(t/s), \vdash \Delta(t/s), \Delta', s\neq t \quad \mbox{iff} \quad\Gamma', \Gamma \vdash \Delta, \Delta'
$$

It is useful to notice the following equivalence
\begin{proposition}
The sequent  $z\in D$ and $(\exists x\in D)z=x$ are equivalent.
\end{proposition}
Proof: derivation of $(\exists x\in D)z=x\vdash z\in D$: $\vdash z\in D, (z\in D)^d$ is equivalent to $z=y\vdash z\in D, (y\in D)^d$, that is $(\exists x\in D)z=x\vdash z\in D$. Derivation of 
$z\in D\vdash (\exists x\in D)z=x$: from the axioms $\vdash z=z$ and $z\in D, (z\in D)^d\vdash $ by the rule $\exists r$ given above. 
\medbreak
One direction of the above equivalence is also propositional:
\begin{proposition}\label{sempre}
The sequent $z=t_1\vee\dots\vee z=t_n\vdash z\in D$ is provable for every domain $D=\{t_1, \dots ,t_n\}$. 
\end{proposition}
Proof: 
The axiom $t_i\in D\vdash t_i\in D$ is equivalent to $t_i\in D, z=t_i\vdash z\in D$ by definition of $=$. Then, cutting the true assumption $t_i\in D$,
one has $z=t_i\vdash z\in D$ for all $i=1\dots n$, namely  $z=t_1\vee\dots\vee z=t_n\vdash z\in D$ by definition of $\vee$.
\medbreak
We say that the domain $D=\{t_1, \dots ,t_n\}$ is {\em focused} when the converse sequent holds, namely when it is
$$
z\in D\vdash z=t_1\vee\dots\vee z=t_n
$$
So we could say that the information \av $z\in D$" has an infinitary content, that becomes finite when $D$ is focused.
 
We introduce the following rule for
the substitution of the free variable $z$ by a term $t$ denoting an element of the domain $D$:
$$
\frac{\Gamma\vdash \Delta}{\Gamma(z/t)\vdash \Delta(z/t)}\,subst(D)
$$
that is labelled by $D$ since, in our model, the validity of substitution depends on $D$.

We see that $subst(D)$ is  derivable, for elements of focused sets, in the following terms:
\begin{lemma}\label{lemmasost}
Let us assume the set $\{t_1,\dots ,t_n\}$ is focused. 
The sequent $\Gamma, z\in D\vdash A(z)$ is equivalent to the $n$ sequents $\Gamma\vdash A(t_i)$, and then to the $n$ sequents $\Gamma\vdash A(t_i)$, by weakening
of the assumption $t_i\in D$
\end{lemma}
$\Gamma, z\in D\vdash A(z)$ means $\Gamma, z=t_1\vee\dots\vee z=t_n\vdash A(z)$, because $D$ is focused. The last is $\Gamma, z=t_i\vdash A(z)$ for every $i$, by definition of $\vee$. These are equivalent to $\Gamma\vdash A(t_i)$ for every $i$ by definition of equality.
\medbreak

Then one can prove the following characterization (see \cite{Ba2}):
\begin{proposition}\label{focuschar}
Let us consider a non empty domain $D=\{t_1, \dots ,t_n\}$.
The sequent $A(t_1)\&\dots \& A(t_n)\vdash (\forall x\in D)A(x)$ is provable for every formula $A$ if and only if $D$ is focused.
\end{proposition}
Proof: 
Let us assume that $D$ is focused.  The sequent $\Gamma\vdash A(t_1)\&\dots \& A(t_n)$ is equivalent to the $n$ sequents $\Gamma\vdash A(t_i)$. By the above lemma, this means $\Gamma, z\in D\vdash A(z)$. For every $\Gamma$ which does not contain $z$ free, this is $\Gamma\vdash (\forall x\in D)A(x)$,
 then in particular for $\Gamma= A(t_1)\&\dots \& A(t_n)$. 

\noindent As for the \av only if": Let us consider $A(x, y) \,=\, x\neq y$. Then, by hypothesis, it is  $z\neq t_1 \&\dots \& z\neq t_m \vdash (\forall x\in D)z\neq x$. This is equivalent to
$z\neq t_1 \&\dots \& z\neq t_m , y\in D\vdash z\neq y$, that, by duality, gives $y\in D, z=y\vdash z=t_1\vee \dots \vee z=t_m$, from which one derives
$(\exists x\in D)z=x\vdash z=t_1\vee \dots \vee z=t_m$. Since  $z\in D\vdash (\exists x\in D)z=x$,
one has $z\in D\vdash z=t_1\vee \dots \vee z=t_m$ cutting the existential formula. 
\medbreak
The converse direction of the sequent is always provable, by substitution:
\begin{proposition}\label{subst}
The sequent $(\forall x\in D)A(x)\vdash A(t_i)$ is provable for every non empty domain $D=\{t_1, \dots ,t_n\}$. Hence $(\forall x\in D)A(x)\vdash A(t_1)\& \dots A(t_n)$ is provable. Symmetrically it is $A(t_1)\vee\dots\vee A(t_n)\vdash (\exists x\in D)A(x)$. Then, the sequent $(\forall x\in D)A(x)\vdash (\exists x\in D)A(x)$ is provable for every non empty domain $D$.
\end{proposition}
Proof:
From $A(z)\vdash A(z)$ and $z\in D\vdash z\in D$ one derives the sequent $(\forall z\in D)A(z), z\in D\vdash A(z)$. By substituting $z$ by $t_i$ in it, and then cutting the true assumption $t_i\in D$, one has $(\forall x\in D)A(x)\vdash A(t_i)$, for every $i=1\dots n$. Then it is $(\forall x\in D)A(x)\vdash A(t_1)\& \dots A(t_n)$ by definition of $\&$. Symmetrically for $\exists$. The last sequent  is proved from $(\forall x\in D)A(x)\vdash A(t_i)$ and its symmetric $A(t_i)\vdash
(\exists x\in D)A(x)$,
by cutting the propositional formula $A(t_i)$.
\medbreak
In particular, in lemma \ref{lemmasost} we have seen that the  substitutions performed in the above proof are valid if $D$ is focused. 
Summing up, we can say that:
\begin{proposition}
$(\forall x\in D)A(x)$ is equivalent to the propositional formula $A(t_1)\& \dots \& A(t_n)$ if and only if the domain $D$ is focused.
\end{proposition}

We assume that singletons are focused, namely
$$
z\in \{u\}\vdash z=u
$$
for every singleton $\{u\}$. Such an assumption is due to our extensional concept of set. Then, $z\in \{u\}$ is a synonimous of $z=u$, since $z=u\vdash z\in \{u\}$ is an instance of \ref{sempre}. 

By the assumption on extensionality of singletons, and the above propositions, it is immediate to see that the metalinguistic link {\em forall} behaves in a symmetric way when the domain is a singleton:
\begin{proposition}\label{singl} 
$A(u)=(\forall x\in \{u\})A(x)$ and $(\exists x\in x\in \{u\})A(x)=A(u)$ are  provable for every $A$, for every singleton $\{u\}$.
So is $(\exists x \in \{u\})A(x)=(\forall x \in \{u\})A(x)$ for every $A$, for every singleton $\{u\}$.  
\end{proposition}
Proof: The first two statements are an immediate consequence of proposition \ref{focuschar}; the third is  immediate by transitivity, that means cutting the proposiional formula $A(u)$. 
It is also very useful to see a direct proof of the sequent $(\exists x\in \{u\})A(x)\vdash (\forall x\in \{u\})A(x)$, without cut. It  is equivalent to the
sequent $A(y), z\in \{u\}\vdash A(z), y\notin \{u\}$, by definition of $\forall$ and $\exists$. The last means 
$$
z=u, A(y)\vdash A(z), y\neq u
$$ 
that in turn is derivable from the axiom $A(u)\vdash A(u)$ by definition of $=$ and $\neq$. 
\medbreak

\subsection{Virtual singletons}\label{virtualsingl}
We extend our notion of singleton, in order to characterize a wider class of symmetric objects. 
Let us term  {\em virtual singleton} any set $V$ for which it
is $(\exists x \in V)A(x)=(\forall x \in V)A(x)$ for every $A$. In particular, we shall term {\em extensional singleton} a set of the form $\{u\}$, namely a set $V$ for which the language has a closed term $u$ such that $z\in V$ iff $z=u$. We have just seen that extensional singletons are virtual singletons,since it is true that
\begin{equation}\label{extax}
z=u, A(y)\vdash A(z), y\neq u
\end{equation}
In order to characterize virtual singletons in general,
\iffalse
The proof of such sequents given above requires substitution.
Let us assume that there are domains $D$ inhabited by some variable, namely that $z\in D$ is true even if the term $z$ is a variable. It will be the case in our model, where a variable is a random variable and $D$ consists of its outcomes. Then $D$ is effective (non empty) in spite of the fact that the outcomes of the random variable are given. Then 
the above result is proved without any occurrence of a substitution rule:
\begin{proposition}\label{dom}
The sequent $(\forall x\in D)A(x)\vdash (\exists x\in D)A(x)$ is provable for every non empty domain $D$.
\end{proposition}
One derives the sequent $z\in D, (\forall x\in D)A(x)\vdash (\exists x\in D)A(x), (z\in D)^d$ by the rules $\forall r$ and $\exists r$ given above, applied to the axioms
$z\in D\vdash z\in D$, $A(z)\vdash A(z)$ and $(z\in D)^d\vdash (z\in D)^d$. Then the formulae $z\in D$ at the left and $(z\in D)^d$ at the right can be cut by our hypothesis on non-emptyness of $D$. 
\medbreak
\fi
we rewrite such sequents, substituting $z=u$, $y\neq u$ by $z\in V$, $(y\in V)^d$, where $d$ is a duality. We have the schemata:
\begin{equation}\label{virtax}
z\in V, A(y)\vdash A(z), (y\in V)^d
\end{equation}
for every $A$.
Such sequents are, in general, underivable, and we shall term them $d${\em -axioms}. One immediately derives $(\exists x \in V)A(x)\vdash(\forall x \in V)A(x)$ from them.
As for the converse direction $(\forall x \in V)A(x)\vdash (\exists x \in V)A(x)$, it is derivable as soon as $V$ is not empty, but its derivation would require substitution. We shall discuss the topic of non emptyness for our model.
Virtual singletons are then non-empty sets $V$ for which the logical system has a suitable duality $d$ for which $d$-axioms can be put.

One proves the following important fact:
\begin{proposition}
Let us assume that a set $V$ has $d$-axioms. If the system is closed by  substitution  by elements of $V$, then the system can prove that $V$ is a singleton.
\end{proposition}
Let $u, u'\in V$ and let us consider $A(z)\equiv z=u$. Applying the substitutions $z/u'$ and $y/u$, one has $u'\in V, u=u\vdash u'=u, (u\in V)^d$.
Since reflexivity $\vdash u=u$ is derivable by the definition of equality, and the sequents $\vdash u'\in V$ and $(u\in V)^d\vdash$ hold by hypothesis, one derives $\vdash u'=u$ by cut.
\medbreak

By the above proposition and lemma \ref{lemmasost}, we have:
\begin{proposition}
 A focused set has $d$-axioms if and only if it is  an extensional singleton. 
 \end{proposition}
 
Then a virtual singleton is either an extensional singleton or an unfocused domain $V$ in a system equipped with a suitable duality $d$ and without substitution rule for elements of $V$.

As soon as a substitution rule of free variables by closed terms denoting elements of $V$ is included in the system, in order to keep consistency, we have two choices:

\noindent 1) we need to drop the duality $d$

\noindent 2)
 we need to consider a new closed term $u$ of the language such that $V=\{u\}$. This identifies $d$ with the duality given by $=$ and $\neq$, the $d$-axioms with the sequents ((\ref{extax}), and the formulae with the propositional formulae $A(u)$. 

On the contrary, virtual singletons, in the generic case, lack of a closed term to describe their content, so we are forced to a predicative description, by a unique quantifier (usually we shall adopt $\forall$). We now see how this can give us an advantage in describing pure symmetry.

\subsection{Extending the action of virtual singletons}\label{entanglement}
The analogy with the behaviour of singletons can be furtherly extended to  virtual singletons. We first observe an immediate consequence of $d$-axioms. Let us consider the sequent $\Gamma(-z, -y), z\in V\vdash A(z), B(z)$, where $V$
is a virtual singleton. Assuming the $d$-axiom $y\in V, A(z)\vdash A(y), (z\in V)^d$, one can cut the formula $A(z)$ and obtain 
$\Gamma(-z, -y), z\in V, y\in V \vdash A(y), B(z), (z\in V)^d$. This is classically equivalent to $\Gamma(-z, -y), z\in V, y\in V \vdash A(y), B(z)$, that is $\Gamma(-z, -y)\vdash (\forall x\in V)A(x), (\forall x\in V)B(x)$, by definition of $\forall$, and then equivalent to $\Gamma(-z, -y)\vdash (\forall x\in V)A(x)\ast (\forall x\in V)B(x)$, interpreting the comma at the right as the multiplicative disjunction, here denoted  $\ast$.
In particular, one proves $(\forall x\in V)A(x)\ast B(x)\vdash (\forall x\in V)A(x)\ast (\forall x\in V)B(x)$, putting $\Gamma=(\forall x\in V)A(x)\ast B(x)$. As is well known, the converse sequent is  derivable for any domain,  then one has the equality
$$
(\forall x\in V)A(x)\ast B(x)=(\forall x\in V)A(x)\ast (\forall x\in V)B(x)
$$
for every virtual singleton $V$, and every pair of formulae $A$, $B$.

Such an equality is not sound in logic, in general. An exception is when $V=\{u\}$ is an extensional  singleton. In such a case, both sides of the equality are equivalent to the propositional compound
formula $A(u)\ast B(u)$, where the propositional connective $\ast$, the multiplicative disjunction, which is not symmetric, appears. The compound objects so described are not symmetric objects.

For virtual singletons, one should introduce a different perspective which avoids inconsistency in the usual logical
calculus, extending symmetry in a purely predicative framework. 

We consider again the pattern of extensional singletons.
For them, $d$-axioms have the form
$z=u, A(y)\vdash A(z), y\neq u$. They are sound because of transitivity of the equality relation: one has that $z=u$ and $y=u$ yield $z=y$. By analogy, let us assume that the premises $z\in V$ and $y\in V$, together, have as a consequence 
a correlation between $z$ and $y$, induced by the virtual singleton $V$. The correlation is an equality when both $z$ and $y$ are equal to a constant $u$. So, if $V$ is a virtual singleton, a  sequent of the form  $\Gamma, z\in V, y\in V\vdash A(z), B(y)$ hides an additional information, that is the correlation between the variables $z$ and $y$.  The comma \av $,$", in a sequent, does not consider the correlation: for, usually, logical rules can be applied to $A$ (resp. $B$), keeping $B$ (resp. $A$) as a context, hence the pieces of information contained in $A$ and $B$ are processed independently. How could logic take care of the correlation between $A$ and $B$ induced induced by that between $z$ and $y$?

\iffalse
In particular, substitution cannot be applied to the variables $z$ and $y$ independently, since they are correlated. 

For our pourposes, we need to consider only the case  it is $y=f(z)$, $f$ an invertible function. Then
we need to consider only one free variable and study  assertions of the following form:
$$
\Gamma(-z), z\in V \vdash A(z),_V B(z)
$$
Substituting the variable $z$ by a closed term $u$ will cause the collapse of the indexed comma into a simple comma:
$$
\frac{\Gamma(-z), z\in V \vdash A(z),_V B(z)}{\Gamma(-z), u\in V \vdash A(u), B(u)}\,subst(V)
$$
Then, in particular, the correlation index is not relevant when $V=\{u\}$ is an extensional singleton, as expected.
\fi

In order to better focus the logical meaning of the correlation and the problem of contexts, we try to link the correlation of first order variables to the correlation of formulae. To this aim, we adopt the notation of indexed formulae: indexes are used as \av second order variables" to distinguish  formulae: $A_1, A_2,\dots A_i,A_j,\dots$. 
Moreover, for our pourposes, we need to  consider  only the case  $y=f(z)$, where $f$ is an invertible function. Then we need to consider only one free variable $z$. So we study  assertions of the following form:
$$
\Gamma(-z), z\in V \vdash A_i(z),_f A_j(z)
$$
where $V$ is a virtual sigleton and where $,_f$ indicates the correlation between the formulae $A_i$ and $A_j$ induced by the premise $z\in V$ by means of the function $f$.
%If $V$ is $\{u\}$, the correlation induced by $f$ is not meaningful, since the domain of $f$ is a singleton, so $,_f$ is simply a comma.

We reconsider the problematic equality, now written with indexes:
$$
(\forall x\in D)A_1(x)\ast A_2(x)=(\forall x\in D)A_1(x)\ast (\forall x\in D)A_2(x)
$$
The equality is true for every domain $D$, when $A_1=A_2=A$, if and only if the disjunction $\ast$ is idempotent, namely $A\ast A=A$ for every $A$. In our perspective, the idea is that the set of indexes
$\{1, 2\}$ is like a virtual singleton, if we cannot distinguish $1$ and $2$. We would like to develop this point. First, one can see that the idempotency of $\ast$ can be derived from
the equivalence:
$$
\Gamma\vdash A, A \qquad \mbox {iff}\qquad \Gamma\vdash A
$$
and conversely.
The  equivalence permits to derive the following, more general:
$$
\Gamma\vdash A(y), A(z) \qquad \mbox {iff}\qquad \Gamma, z=y\vdash A(z)
$$
since both sides are equivalent to $\Gamma, z=y\vdash A(z), A(z)$, by definition of $=$ and by idempotency, respectively. In particular, one has exactly idempotency, considering twice a variable $z$ which does not appear in $A$, since $z=z$ is true.
  
We now lift the last equivalence to the second order, namely to indexes $i,j$. We write
$$
\Gamma\vdash A_i, A_j \qquad \mbox {iff}\qquad \Gamma, i=j\vdash A_i
$$
Again this is idempotency when $i$ and $j$ coincide, given that $i=i$ is true. 

In our perspective, we generalize the equivalence to virtual singletons, as follows:
$$
\Gamma\vdash A_i,_f A_j \qquad \mbox {iff}\qquad \Gamma, i\sim_f j\vdash A_i
$$
where $A_i=A_i(z)$, $A_j=A_j(z)$, $\Gamma=\Gamma', z\in V$. We term such an equivalence {\em second order conversion}.
The correlation $,_f$ given by the virtual singleton $V$ and by $f$, is translated into a relation $\sim_f $ between the indexes $i$ and $j$. We are saying that we cannot distinguish two indexes (in particular, they could be simply the same!) This extends virtual singletons to sets of indexes $I$: $i\sim_f j$ is $i\in I$ and $j\in I$, where $I$ is a virtual singleton of indexes.
For every $V$ and $f$, we put the definitory equation of a connective $\Join_f$ interpreting the correlation $,_f$ between formulae. It imitates the definitory equation of $\ast$:
$$
\Gamma', z\in V\vdash A_i \Join_f A_j \quad\equiv\quad \Gamma', z\in V\vdash A_i,_f A_j
$$ 
But then it is 
$$
\Gamma', z\in V\vdash A_i \Join_f A_j \quad\mbox{iff} \quad \Gamma', z\in V, i\sim_f j\vdash A_i
$$
by second order conversion. So $\Join$ is also like a quantifier, whose domain is an unfocused set of indexes $I$, where $i, j\in I$.

We now see that $\Join$  commutes w.r.t the quantifier $\forall$, generalizing the case of extensional singletons in this way. 
We make the assumption that indexing a formula is independent of its logical construction. This means that if any connective is applied to an indexed formula $A_i$, the result is still indexed by $i$ and conversely. Then one could prove that the construction of formulae by logical rules occurs \av in parallel" rather than in a context-free construction. In particular, we are interested in parallelizing the formation rule of $\forall$ (see \cite{Ba}): 
\iffalse
\begin{proposition}
Let us assume to have a connective $\cal C$ defined putting an equivalence $\Gamma\vdash {\cal C} A \,\mbox{iff}\, \Gamma'\vdash A$. 
If the rule of second order conversion holds, for every virtual singleton $V$, the $\cal C$
rules are applied in a parallel way to $A_i$ and $A_j$ in $\Gamma\vdash A_i,_f A_j$, namely it is $\Gamma'\vdash A_i,_f A_j$ if and only if $\Gamma\vdash {\cal C} A_i,_f {\cal C} A_j$, for every connective $\cal C$.
\end{proposition}
Proof: $\Gamma'\vdash A_i,_f A_j$ is equivalent to $\Gamma', i\sim j\vdash A_i$ by second order conversion. This is equivalent to $\Gamma, i\sim j\vdash {\cal C} A_i$ by definition of $\cal C$. By our assumption on indexing, ${\cal C} A_i$ is $[{\cal C} A]_i$. Again by conversion, we have equivalently 
$\Gamma\vdash [{\cal C} A]_i,_f [{\cal C} A]_j$, that is $\Gamma\vdash {\cal C} A_i,_f {\cal C} A_j$.
\medbreak
\fi
\begin{proposition}
The parallel rule of $\forall$ with respect to the indexed comma $,_f$:
$$
\frac{\Gamma, z\in V\vdash A_i(z),_f A_j(z)}{\Gamma\vdash (\forall x\in V)A_i(x),_f (\forall x\in V)A_j(x)}
$$
is derivable for every virtual singleton $V$.
\end{proposition}
Proof: $\Gamma, z\in V\vdash A_i(z),_f A_j(z)$ is equivalent to $\Gamma', i\sim j\vdash A_i$ by second order conversion. This is equivalent to $\Gamma, i\sim j\vdash (\forall x\in V)A_i(x)$ by definition of $\forall$. By our assumption on indexing, $(\forall x\in V)A_i(x)$ is $[(\forall x\in V)A(x)]_i$. Again by conversion, we have equivalently 
$\Gamma\vdash [(\forall x\in V)A(x)]_i,_f [(\forall x\in V)A(x)]_j$, that is $\Gamma\vdash (\forall x\in V)A_i(x),_f (\forall x\in V)A_j(x)$.
\medbreak
An analogous reasoning could be applied to other connectives: rules for $,_f$ are parallel.

The above proposition allows to derive the sequent
$$
(\forall x\in V)(A_i\Join A_j)(x)\vdash (\forall x\in V)A_i(x)\Join (\forall x\in V)A_j(x)
$$
applying the $\forall$ parallel rule and the other definitory equations when $\Gamma= (\forall x\in V)(A_i\Join A_j)(x)$ 
One could see that the converse sequent is also derivable, adopting for $\Join$ rules analogous to those of $\ast$. But, perhaps, the best way to see this is that quantifiers on virtual singletons are symmetric, then the converse sequent is derived in the symmetric way. Then one derives the equality,
generalizing the case of extensional singletons.

Summing up, we have a predicative symmetric object, defined by means of $\forall$ and $\Join$, as $(\forall x\in V)A_i(x)\Join (\forall x\in V)A_j(x)$ or equivalently as $(\forall x\in V)(A_i\Join A_j)(x)$, which has one virtual singleton as a domain  for first order variables and another for indexes.
This enriches the set of symmetric objects. In particular when $I=\{1\}$ is an extensional singleton, the quantifier is $\forall$ (or equivalently $\exists$). When $D=\{u\}$ is an extensional singleton, we have the propositional case: $\Join$ is $\ast$ (no correlation), and we have no symmetry, as already observed. 

\section{Duality and symmetry representing qubits}
\subsection{Representation of quantum states}
We briefly remind the predicative representation of quantum states introduced in \cite{Ba}, \cite{Ba2}.
Let us consider a  random variable $Z$, with outcomes $s_1,\dots ,s_m$ and frequencies $p\{Z=s_1\},\dots ,p\{Z=s_m\}$. 
It yields a set
$$
D_Z\equiv \{z=(s(z), p\{Z=s(z)\})\}
$$
where $s(z)$ is the generic outcome and $p\{Z=s(z)\}> 0$ is its frequency.
We term $D_Z$ {\em random first order domain} \cite{Ba2}. $D_Z$ is focused if and only if the set of the outcomes of $Z$ is focused,
for it is $z=t_1 \vee \dots \vee z=t_m$ if and only if it is $s(z)=s(t_1)\vee\dots\vee s(z)=s(t_m)$.  
When $D_Z$ is unfocused,  the open predicate $z\in D_Z$, in the free variable $z$,
is a primitive entity, related to the random variable, which cannot be described by closed terms. 

As is well known, a quantum measurement of a particle is a random variable and hence one has a random first order domain associated to it.
Let us fix a particle $\cal A$ and consider a discrete observable whose measurement produces  a finite domain $D_Z=\{t_1,\dots ,t_m\}$, where $t_i=(s(t_i), p\{Z=s(t_i)\})$.
One can write  $A(t_i)$ for the proposition 
\av The particle $\cal A$ is found in state $s(t_i)$ with probability $p\{Z=s(t_i)\}$". Let us  
summarize all the hypothesis concerning the preparation of the quantum state and its measurement into a set of premises $\Gamma$. One has that $\Gamma$ {\em yield} $A(t_i)$, for $i=1\dots m$. Such  $m$ assertions are written $\Gamma\vdash A(t_i)$ as  sequents \cite{SBF}. Then one has equivalently
$\Gamma\vdash A(t_1)\&\dots \& A(t_m)$, where $\&$ is the additive conjunction. The proposition $A(t_1)\&\dots \& A(t_m)$ represents our knowledge of the state after measurement,
namely the probability distribution of the outcomes.

In order to describe the quantum state prior to measurement, one should drop the identification of the states, possible only after measurement. In our terms,
this means to drop the equality that renders $D_Z$ focused. The we need to describe our knowledge by a free variable. We consider the 
proposition $A(z)$: 
\av The particle is in state $s(z)$ with probability $p\{Z=s(z)\}$" for all $z\in D_Z$. If the measurement hypothesis are denoted by $\Gamma$,
(where $\Gamma=\Gamma(-z)$ does not depend on the variable $z$, since the measurement hypothesis cannot depend on its eventual outcome),
it is $\Gamma(-z)$ {\em  yield}  $A(z)$ {\em forall} $z\in D_Z$, namely \av {\em forall} $z\in D_Z$, $\Gamma(-z)\vdash A(z)$". 
So, considering the definitory equation of the universal quantifier (\ref{foralldef}), 
we attribute a state to the particle by 
the proposition 
$$
(\forall x \in D_Z)A(x)
$$

As seen in prop. \ref{subst}, one derives the sequent $(\forall x \in D_Z)A(x)\vdash A(t_1)\&\dots\& A(t_m)$
by substitution. 
In our terms, the sequent says  that the probability distribution follows from the state. In quantum mechanics, the probability distribution is derived by measurement.
Hence
a measurement is represented by a substitution of the free variable $z$ by the closed terms $t_i$ in our model. 

The converse sequent $A(t_1)\&\dots\& A(t_m) \vdash  (\forall x \in D_Z)A(x)$ 
holds
if and only if $D_Z$ is focused (see \ref{focuschar}). This enables us to characterize quantum states predicatively. 
But, in particular, when the measurement has a unique certain outcome (sharp state), namely the random first order domain of the state of the particle is a singleton, its representation is also propositional.

\subsection{Qubits and duality}
We now consider the items of information contained in a quantum particle with respect to a two-valued observable.
We consider the measurement of  the spin of a particle  w.r.t. a fixed axis, say the $z$ axis.  
The outcome of a measurement of a particle $q$  
is \av spin down"  with probability $\alpha^2$ and \av spin up"  with probability $\beta^2$, $\alpha^2+\beta^2=1$.

We consider the usual representation of qubits as vectors in the Hilbert space $C^2$, up to a global phase factor. We fix the orthonormal basis $\{\ket \downarrow, \ket \uparrow \}$. Then the state of a  qubit $q$ is written 
$$
\ket q= \alpha\ket\downarrow + \beta e^{i\phi}\ket\uparrow 
$$ 
Different qubits yielding the same probability distribution can be characterized by a phase $\phi$. The quantum measurement ignores $\phi$\footnote{in order to  consider $\phi$ as a real item of information, one should apply a phase estimation procedure.} and characterizes the real probabilities, given by $\alpha^2$ and $\beta^2$.

The random first order domain of the state represented by $\ket q$ is the set $D_Z=\{(\downarrow, \alpha^2), (\uparrow, \beta^2)\}$. Then the state of the qubit
is represented by the predicative formula
$(\forall x\in \{(\downarrow, \alpha^2), (\uparrow, \beta^2)\})A(x)$. 
The unfocused domain $D_Z=\{(\downarrow, \alpha^2), (\uparrow, \beta^2)\}$,  corresponds to a family
of vectors $\alpha\ket \downarrow + e^{i\phi}\beta\ket \uparrow$, $\phi\in [0, 2\pi)$.  Let us consider two qubits in the same family, $\ket q= \alpha\ket \downarrow+ e^{i\phi}\beta\ket \uparrow$ and $\ket {q'}= \alpha\ket \downarrow+ e^{i\phi'}\beta\ket \uparrow$. We would like to characterize them in our setting
of fixed measurement w.r.t. the basis $\ket \downarrow$, $\ket \uparrow$. It is well known that any two qubits can be distinguished by measurement if and only if they are  orthogonal (see \cite{NC}). The inner product of $\ket q$ and $\ket {q'}$ in $C^2$ is $\alpha^2+ e^{i(\phi'-\phi)}\beta^2$. It is $\alpha^2+ e^{i(\phi'-\phi)}\beta^2=0$ if and only if $\alpha^2=\beta^2=1/2$  and $\phi'-\phi=\pi$. We consider the phases $\phi=0$ and $\phi=\pi$, which give real factors, and characterize the couple of orthogonal vectors $\ket +=1/\sqrt 2\ket \downarrow+ 1/\sqrt 2\ket \uparrow$ and $\ket -= 1/\sqrt 2\ket \downarrow- 1/\sqrt 2\ket \uparrow$.

The elements of the  dual basis $\ket +, \ket -$ so obtained, are characterized as the eigenvectors of the Pauli matrix 
$X= \left (
\begin{array}{cc}
0 & 1 \\
1 & 0
\end{array}
\right )$ which corresponds to the observable \av spin along the $x$ axis" in our hypothesis. It is $X\ket +=\ket +$ and $X\ket -=\ket -$ (up to a global phase factor). The elements of the computational basis  are eigenvectors of the Pauli matrix 
$Z= \left ( \begin{array}{cc}
1 & 0 \\
0 & -1
\end{array}
\right )$ that is the observable \av spin along the $z$ axis: it is $Z\ket\downarrow=\ket \downarrow$ and $Z\ket \uparrow=\ket\uparrow$ (up to a global phase factor). Moreover, $X$ switches the elements of the given computational basis for our measurements: $\downarrow$ and $\uparrow$:
$$
X\ket \downarrow =\uparrow \qquad X\ket \uparrow =\downarrow
$$
Hence $X$ is the $NOT$ gate for our computation. 
On the contrary, $Z$ switches the elements of the dual basis, namely the phases:
$$
Z\ket +=\ket - \qquad Z\ket - = \ket +
$$
 
So a complex duality is definable on quantum states, however only one half of it can emerge by measurement: for, the observables \av spin along $z$" and \av spin along $x$" are incompatible.
The action of $Z$ is non trivial only on the dual basis. Then we can see only the duality given by $NOT$ after measurement. 

\subsection{Qubits as virtual singletons}
We describe the duality given by the $NOT$ and by the $Z$ gate 
in logic. 
The measurement of qubits in the computational basis $\ket\downarrow, \ket\uparrow$ yields the couple of extensional singletons 
$$
D_\downarrow=\{(\ket \downarrow, 1)\} \qquad   D_\uparrow=\{(\ket \uparrow, 1)\}
$$

The measurement of particles in the dual basis $\ket +, \ket -$
yields the domain $D=\{(\downarrow, 1/2), (\uparrow, 1/2)\}$, given by the uniform distribution of the outcomes. In our measurement context, particles whose state is such cannot be given an objective property. Their characterization depends on the phase $+$ or $-$. So we have two unfocused copies of $D$: 
$$
D_+ \qquad D_-
$$ 
We assume that both are non empty, even if they are not accessible by substitution, since they correspond to a quantum state, that would
be described by a singleton assuming a different measurement context. We think that they are inhabited by the uniform random variable. 
In our measurement context, $D_+$ and $D_-$ are conceivable as virtual singletons. They are equal as sets, from an extensional point of view. The labels
$+$ and $-$ are like two modalities which can characterize the attribution of one of the two orthogonal states $\ket +$ and $\ket -$ to a certain qubit.

We transfer the correspondences between quantum states given by $NOT$ and $Z$ into logic, putting the following dualities  $\perp$  and $\top$, respectively, on the above domains:
$$
D_\downarrow^\perp\equiv D_\uparrow \quad D_\uparrow^\perp\equiv D_\downarrow
\qquad
D_\downarrow^\top \equiv D_\downarrow \quad D_\uparrow^\top\equiv D_\uparrow 
$$
$$
D_+^\top \equiv D_- \quad D_-^\top \equiv D_+
\qquad
D_+^\perp \equiv D_+ \quad D_-^\perp \equiv D_-
$$
It is $z\in D_\downarrow$ iff $z=\downarrow$, and also $(z\in D_\downarrow^\perp\equiv z\in D_\uparrow$ iff $z\neq \downarrow$, since any two quantum states can be distinguished if and only if they are orthogonal. Then $z\in D_\downarrow^\perp$ is a dual proposition w.r.t $z\in D_\downarrow$.
This pattern can be extended to the duality $\top$: the dual of $z\in D_+$ is $z\in D_-$ and conversely. So we put the following  $\top$-axioms:
$$
z\in V, A(y)\vdash A(z), y\in V^\top
$$
for $V=D_+$ and $V=D_-$.

The extensional singletons $D_\downarrow$ and $D_\uparrow$ form the propositions $(\forall x\in D_\downarrow)A(x)$ and $(\forall x\in D_\uparrow)A(x)$. They are equivalent to the
propositional formulae $A(\ket \downarrow, 1)$ and $A(\ket \uparrow, 1)$, abbreviated $A_\downarrow$ and $A_\uparrow$, respectively. 

The virtual singletons $D^+$ and $D^-$ form the propositions $(\forall x\in D_+)A(x)$ and $(\forall x\in D_-)A(x)$, abbreviated $A_+$ and $A_-$ respectively.

In our model, such propositions attribute a state to a particle. For different particles $\cal A, \cal B, \ldots$, we have different propositions. Let us consider two lists of couples of literals of this kind: 
$$
A_\downarrow, A_\uparrow; B_\downarrow, B_\uparrow, \dots
$$ 
({\em sharp literals}); and 
$$
A_+, A_-; B_+, B_-, \dots
$$ 
({\em phase literals}). Sharp literals are recognizable as couples of opposites, putting the duality $(-)^\perp$ on them as follows:
$$
A_\downarrow^\perp\equiv A_\uparrow \qquad\qquad A_\uparrow^\perp \equiv A_\downarrow
$$
We can prove that
\begin{proposition}
The definition of $\perp$  is compatible with a negation in the predicative logical language.
\end{proposition}
We have: $[(\forall x\in D_\downarrow)A(x)]^\perp =A_\downarrow^\perp = A_\uparrow= (\forall x\in D_\uparrow)A(x)$, so our position implies $[(\forall x\in D_\downarrow)A(x)]^\perp = 
(\forall x\in D_\downarrow^\perp)A(x)$. On the other side, let us consider the usual classical definition of negation for quantified formulae: $[({\cal Q} x\in D)A(x)]^\perp
=({\cal Q}^\perp x\in D)A(x)^\perp$ for any quantifier $\cal Q$ and any domain $D$. Moreover, let us put: $({\cal Q} x\in D_\downarrow)A(x)^\perp\equiv ({\cal Q} x\in D_\uparrow)A(x)$. Namely, we assume that $A(x)^\perp$ is $A(x)$, where $x$ can be in one of the two 
opposite domains $D_\downarrow$ and $D_\uparrow$.
 One has that $[(\forall x\in D)A(x)]^\perp$ is $(\forall x\in D^\perp)A(x)$ by our definition. On the other side, it is  $[(\forall x\in D)A(x)]^\perp=(\exists x\in D)A(x)^\perp=
 (\exists x\in D^\perp)A(x)$ by the above assumption. Then, since it is $(\exists x\in D^\perp)A(x)=(\forall x\in D^\perp)A(x)$, one can conclude. 
\medbreak 
We put the analogous definition for phase literals and the duality $\top$:
$$
(A_+)^\top \equiv A_- \quad (A_-)^\top \equiv A_+
$$

Extending the whole duality given by the $Z$ and the $NOT$ quantum gates, described above, to literals, we put also
$$
A_\downarrow^\top\equiv A_\downarrow \quad A_\uparrow^\top\equiv A_\uparrow 
\qquad\qquad 
((A_+)^\perp\equiv A_+ \quad (A_-)^\perp\equiv A_-
$$
Phase literals are our symmetric literals: they are fixed points for the duality $\perp$.  
  
The observations at the end of section \ref{virtualsingl} can now be applied to our model. 
Measuring a particle $\cal A$ in state $\ket +$ or $\ket -$ in our original measurement context, determines its collapse into the mixed state described by the propositional formula $A(\downarrow, 1/2)\& A(\uparrow, 1/2)$, in both cases. This determines the loss of the duality $\top$. On the contrary, switching the measurement context
and considering the observable \av spin w.r.t. the $x$ axis" gives an objective property to the state of $\cal A$. This can be described by the identification
of the domains $D_+$ and $D_-$ with two extensional singletons, let us label them $\{+\}$  and $\{-\}$, where $+$ and $-$ are new closed terms to denote the state
to be attributed to the particle. Then the duality $\top$ is identified with that given by the equalities $z=+/z=-$, and the duality $\perp$ disappears.

\subsection{Extending phase literals and symmetry}
We now consider a couple of particles ${\cal A}_1, {\cal A}_2$ in one of the four Bell's states:
$$
1/\sqrt 2\ket{\downarrow\downarrow}\pm 1/\sqrt 2\ket{\uparrow\uparrow}\qquad 1/\sqrt 2\ket{\downarrow\uparrow}\pm 1/\sqrt 2\ket{\uparrow\downarrow}
$$
For each state, one finds the state of the second particle correlated to the outcome of the measurement of the first: it can be identical or opposite. The correlation
is due to a correlation between the two particles. 

As seen in Ghirardi (\cite{Ghi}, p. 306), it is not limiting to assume that the entangled state described  is induced by two identical particles, since such an assumption leads, formally,  to the Bell's states in a natural way. So it is natural for us to require that, if $j\in I=\{1,2\}$ is the index of one of the two, one cannot say which one it is. So we adopt the formalism we have studied in subsection \ref{entanglement} to represent the four Bell's states. 
We write down the following four cases for the state of the two particles:
$$
\Gamma, z\in D_+\vdash A_1(z),_i A_2(z) \qquad \Gamma, z\in D_-\vdash A_1(z),_i A_2(z)
$$
$$
\Gamma, z\in D_+\vdash A_1(z),_o A_2(z) \qquad \Gamma, z\in D_-\vdash A_1(z),_o A_2(z)
$$
where we have two modalities for the phase: $+$ and $-$, and
where the indexed commas $,_{i/o}$ describe the identical or opposite correlation between the outcomes of measurement for the couple of states. As seen in subsection \ref{entanglement}, the correlations $i$ and $o$ between the outcomes correspond to  the correlations, that we could label again by $i$ ad $o$, on the set of indexes $I=\{1,2\}$. Then we have two modalities for the set $I=\{1,2\}$ indexes too. 

\noindent Then one represents the four Bell's states adopting the generalized quantifier obtained from $\forall$ and $\Join$, with domains $V=D_+, D_-$ for the first order variables and correlations $\Join_i$ and $\Join_o$ for the indexes.
$$
(\forall x\in D_+)(A_1(x)\Join_i A_2(x)) \qquad (\forall x\in D_+)(A_1(x)\Join_o A_2(x))
$$
$$
(\forall x\in D_-)(A_1(x)\Join_i A_2(x)) \qquad (\forall x\in D_-)(A_1(x)\Join_o A_2(x))
$$
The propositions so obtained are a generalization of phase literals, since they consist of a generalized symmetric quantifier applied to virtual singletons.

We briefly discuss how to extend the dualities  $\perp$ and $\top$ to the new propositions. 
Let us consider $\perp$. One should keep that the right way to extend $\perp$ to them is the identity, since they are like phase literals. Moreover, the duality $\perp$ is the identity on the domains and on the connectives $\forall$ and $\Join$, since they are symmetric. Then, even applying compositionally the duality to our propositions, one finds the identity, that  extends the identity induced by $\perp$ on phase literals. 

On the other side, let us consider $\top$. What is the right way to extend it to our propositions? Applying compositionally the duality $\top$ to them, one would not find the identity, since $\top$ switches $+$ and $-$. But one has to consider the domain of indexes too. If again $\top$ induces a switch, one should exchange the two modalities $i$ and $o$. This does not correspond to the real facts. For example, the singlet state 
$1/\sqrt 2\ket{\downarrow\uparrow} - 1/\sqrt 2\ket{\uparrow\downarrow}$ is a fixed point for the spin  with respect to any axis. Then we think that the compositional application is not a good attitude in order to find the right way to extend $\top$. Actually, our propositions are like literals and have no proper subformula. We could simply keep that, since the observable $Z$ induces a switch on modalities, the global effect on the two domains, of variables and of indexes,  together, is to have no switch. Then $\top$ would be the identity on these formulae. This means that no change of measurement context is able to let a negation emerge in this case. Bell's states are really \av hidden" logical objects, and their hidden presence is really against the direction of logical consequence, whatever measurement context one assumes.


\begin{thebibliography}{99}
\bibitem[Ba]{Ba} 
Battilotti, G., Interpreting quantum parallelism by sequents, International Journal of Theoretical Physics 49 (2010) 3022-3029.
\bibitem[Ba2]{Ba2} 
Battilotti, G., Characterization of quantum states in predicative logic, International Journal of Theoretical Physics 50 (2011) 3669-3681.
\bibitem[BS]{BS}
Baltag, A. Smets, S. Correlated Knowledge: an epistemic-logic view of quantum entanglement,
International Journal of Theoretical Physics, Online July 2010.
\bibitem[BV]{BV}
Bonzio, S., Verrucchi, P., Open Quantum Systems dynamics and quantum algorithms, submitted.
\bibitem[Ca]{Ca} 
Castagnoli, G. Probing the mechanism of quantum speed up by time-symmetric quantum mechanics, arXiv quant-ph/1107.0934v7.
\bibitem[DCGL03]{DCGL03}
Dalla Chiara M.L.,  Giuntini R.,  Leporini R., Quantum computational logics. A survey,  in V. F. Hendricks, J. Malinowski eds., Trends in Logic: 50 years
of studia logica, Kluwer Academic Publishers, Dordrecht (2003) 213-255.
\bibitem[DCGL]{DCGL2}
Dalla Chiara M.L.,  Giuntini R.,  Leporini R., 
Compositional and Holistic Quantum Computational Semantics,  Natural Computing 6  (2007)  113-132.  
\bibitem[DE]{DE}
Dolev, S., Elitzur A.C., Non sequential behavior of the wave function.  arXiv quant-ph/ 0102109.
\bibitem[FS]{FS} 
Faggian, C., Sambin, G., From basic logic to quantum logic with cut-elimination, Proc. IQSA96, International Journal
of theoretical Physics  31 (1998) 31-37.
\bibitem[Ghi]{Ghi} 
Ghirardi, G.C., {\em Un'occhiata alle carte di Dio}, Il Saggiatore, Milano (2009). 
\bibitem [Ma]{Ma}
Maietti, M.E., Lecture notes in logic, course of Logic for Computer Science, University of Padua.
\bibitem[MS]{MS}
Maietti, M.E., Sambin, G.
Toward a minimalist foundation for constructive mathematics, in ``From Sets and Types to Topology and Analysis: 
Towards Practicable Foundations for Constructive Mathematics" (L. Crosilla, P. Schuster, eds.), Oxford UP, 2005.
\bibitem[MB]{MB}
Matte Blanco, I., {\em The unconscious as infinite sets}, Duckworth, London, 1975. 
\bibitem[NC]{NC}
Nielsen, N. A., Chuang, I. L. {\em Quantum Computation and Quantum Information}, 
Cambridge University Press (2000).
\bibitem[SBF]{SBF}
Sambin G., Battilotti G., Faggian C., Basic logic: reflection, symmetry, visibility,  The Journal of Symbolic Logic 65 (2000) 979-1013.
\end{thebibliography}
\end{document}